\documentclass[11pt]{amsart}

\usepackage{latexsym,amsmath,amsthm,amssymb,amscd,amsfonts,diagrams,graphics}

\DeclareFontFamily{U}{rsf}{}
\DeclareFontShape{U}{rsf}{m}{n}{
  <5> <6> rsfs5 <7> <8> <9> rsfs7 <10-> rsfs10}{}
\DeclareMathAlphabet{\mathscr}{U}{rsf}{m}{n}

\DeclareMathAlphabet{\mathgth}{U}{euf}{m}{n}

\DeclareFontFamily{U}{cyr}{}
\DeclareFontShape{U}{cyr}{m}{n}{
  <5> wncyr5 <6> wncyr6 <7> wncyr7 <8> wncyr8 <9> wncyr9 <10-> wncyr10}{}
\DeclareMathAlphabet{\mathcyr}{U}{cyr}{m}{n}

\input cyracc.def

\newcommand{\cB}{{\mathscr B}}

\newcommand{\sX}{{\mathgth X}}
\newcommand{\cC}{{\mathscr C}}
\newcommand{\cE}{{\mathscr E}}
\newcommand{\cF}{{\mathscr F}}
\newcommand{\cG}{{\mathscr G}}

\newcommand{\cL}{{\mathscr L}}

\newcommand{\cO}{{\mathscr O}}

\newcommand{\hcB}{\widehat{\cB\,\,}\!\!}
\newcommand{\hcC}{\widehat{\cC}}

\newcommand{\FMYX}{\Phi_{Y\ra X}}
\newcommand{\FMXX}{\Phi_{X\ra X}}
\newcommand{\FMXY}{\Phi_{X\ra Y}}

\newcommand{\FMYZ}{\Phi_{Y\ra Z}}
\newcommand{\FMXZ}{\Phi_{X\ra Z}}
\newcommand{\fmYX}{\phi_{Y\ra X}}
\newcommand{\fmXY}{\phi_{X\ra Y}}

\newcommand{\fmYZ}{\phi_{Y\ra Z}}
\newcommand{\fmXZ}{\phi_{X\ra Z}}

\newcommand{\even}{{\mathrm{even}}}

\newcommand{\HKR}{{\mathrm{HKR}}}

\newcommand{\D}{{\mathbf D}_{\mathrm{coh}}^b}

\newcommand{\chk}{{\scriptscriptstyle\vee}}
\newcommand{\R}{\mathbf{R}}
\newcommand{\Ld}{\mathbf{L}}

\DeclareMathOperator{\Tr}{Tr}

\DeclareMathOperator{\Spec}{Spec}

\newcommand{\Hom}{{\mathrm{Hom}}}

\DeclareMathOperator{\Td}{td}
\DeclareMathOperator{\ch}{ch}

\DeclareMathOperator{\Id}{Id}

\DeclareMathOperator{\Ext}{Ext}

\newcommand{\adjoint}{\dashv}

\newcommand{\ra}{\rightarrow}
\newcommand{\lra}{\longrightarrow}

\newcommand{\scdot}{\,\cdot\,}

\newcommand{\C}{\mathbf{C}}

\newcommand{\Q}{\mathbf{Q}}

\newcommand{\gMod}{\mathgth{Mod}}

\newcommand{\iso}{\cong}

\newarrow{Equal} =====

\theoremstyle{plain}
\newtheorem{theorem}{Theorem}[section]

\newtheorem{corollary}[theorem]{Corollary}
\newtheorem{proposition}[theorem]{Proposition}
\newtheorem{conjecture}[theorem]{Conjecture}
\theoremstyle{definition}
\newtheorem{definition}[theorem]{Definition}
\newtheorem{definition-theorem}[theorem]{Definition-Theorem}

\theoremstyle{remark}
\newtheorem{remark}[theorem]{Remark}

\renewcommand{\phi}{\varphi}

\begin{document}

\author{Andrei C\u ald\u araru}

\title[The Mukai pairing, II]{The Mukai pairing, II: the Hochschild-Kostant-Rosenberg isomorphism}

\date{}

\begin{abstract}
We continue the study of the Hochschild structure of a smooth space
that we began in~\cite{CalHH1}, examining implications of the
Hochschild-Kostant-Rosenberg theorem.  The main contributions of the
present paper are:
\begin{enumerate}
\item[--] we introduce a generalization of the usual Mukai pairing on
differential forms that applies to arbitrary manifolds;
\item[--] we give a proof of the fact that the natural Chern character 
map $K_0(X) \ra HH_0(X)$ becomes, after the HKR isomorphism, the usual
one $K_0(X) \ra \bigoplus H^i(X, \Omega_X^i)$; and
\item[--] we present a conjecture that relates the Hochschild and
harmonic structures of a smooth space.
\end{enumerate}
\end{abstract}

\maketitle

\section{Introduction}
\label{sec:intro}

\subsection{}
In~\cite{CalHH1} we introduced the Hochschild structure $(HH^*(X),
HH_*(X))$ of a smooth space $X$, which consists of:
\begin{itemize}
\item[--] a graded ring $HH^*(X)$, the Hochschild cohomology ring,
  defined as
\[ HH^i(X) = \Hom_{\D(X\times X)}(\cO_\Delta, \cO_\Delta[i]), \]
  where $\cO_\Delta=\Delta_*\cO_X$ is the structure sheaf of the
  diagonal in $X\times X$;
\item[--] a graded left $HH^*(X)$-module $HH_*(X)$, the Hochschild
  homology module, defined as
\[ HH_i(X) = \Hom_{\D(X\times X)}(\Delta_!\cO_X[i], \cO_\Delta), \]
  where $\Delta_!$ is the left adjoint of $\Delta^*$ defined by
  Grothendieck-Serre duality~(\cite[3.3]{CalHH1});
\item[--] a non-degenerate pairing $\langle\scdot,\scdot\rangle$
  defined on $HH_*(X)$, the {\em generalized Mukai pairing} (for the
  definition see~\cite{CalHH1}).
\end{itemize}

\subsection{}
\label{subsec:chernchar}
Following ideas of Markarian~\cite{Mar} we also introduced the Chern
character map
\[ \ch:K_0(X) \ra HH_0(X) \]
by setting $\ch(\cF)$ for $\cF\in\D(X)$ to be the unique element of
$HH_0(X)$ such that
\[ \Tr_{X\times X}(\mu\circ \ch(\cF)) = \Tr_X(\FMXX^\mu(\cF)) =
\Tr_X(\pi_{2,*}(\pi_1^*\cF\otimes \mu)) \] 
for every $\mu\in\Hom_{\D(X\times X)}(\cO_\Delta, S_\Delta)$.  

Here $\Tr$ is the Serre duality trace~(\cite[2.3]{CalHH1}), $S_X
= \omega_X[\dim X]$ is the dualizing object of $\D(X)$ (also to be
thought of as the functor $-\,\otimes_X S_X$), $S_\Delta=\Delta_*S_X$
is the object whose associated integral transform is $S_X$, and
$\FMXX^\mu$ is the natural transformation $1_X\Rightarrow S_X$
associated to $\mu$~(\ref{subsec:nattrans}).

It is worth pointing out that $\mu\circ \ch(\cF)$ is a morphism
$\Delta_!\cO_X \ra S_\Delta$, so using the definition of
$\Delta_! = S_{X\times X}^{-1} \Delta_* S_X$ it follows that
$\mu\circ\ch(\cF)$ is in fact a morphism
\[ S_{X\times X}^{-1} S_\Delta \ra S_\Delta, \]
and thus it makes sense to take its trace on $X\times X$.  For more
details see~\cite{CalHH1}.

\subsection{}
The Hochschild structure satisfies the following properties~(\cite{CalHH1}):
\begin{itemize}
\item[--] to every integral functor $\Phi:\D(X)\ra \D(Y)$ there is a
  naturally associated map of graded vector spaces $\Phi_*:HH_*(X)\ra
  HH_*(Y)$.  This association is functorial, commutes with $\ch$, and
  if $\Psi$ is a left adjoint to $\Phi$, then $\Psi_*$ is a left
  adjoint to $\Phi_*$ with respect to the Mukai pairings on $X$ and on
  $Y$, respectively, i.e.,
\[ \langle v, \Phi_* w\rangle_Y = \langle \Psi_*v, w
  \rangle_X \] 
  for $v\in HH_*(Y)$, $w\in HH_*(X)$, and a similar statement holds
  for right adjoints;
\item[--] the Mukai pairing is a generalization of the Euler pairing
  on $K_0(X)$,
\[ \langle \ch(\cE), \ch(\cF) \rangle = \chi(\cE, \cF) = \sum_i (-1)^i
  \dim \Ext^i_X(\cE, \cF) \]
  for any $\cE, \cF\in \D(X)$;
\item[--] the Hochschild structure is invariant under derived
  equivalences given by Fourier-Mukai transforms; in other words, if
  $\FMXY:\D(X)\ra \D(Y)$ is a Fourier-Mukai transform, then there are
  induced isomorphisms $HH^*(X) \iso HH^*(Y)$ (as graded rings),
  $HH_*(X)\iso HH_*(Y)$ (as graded modules over the corresponding
  cohomology rings) and this isomorphism is an isometry with respect
  to the generalized Mukai pairings on $X$ and on $Y$, respectively.
\end{itemize}

\subsection{}
The purpose of this paper is to study the similarities between the
Hochschild structure and the {\em harmonic structure} $(HT^*(X),
H\Omega_*(X))$ of $X$, whose vector space structure is defined as
\begin{align*}
HT^i(X) & = \bigoplus_{p+q=i} H^p(X, \bigwedge^q T_X) \\
H\Omega_i(X) & = \bigoplus_{q-p=i} H^p(X, \Omega^q_X).
\end{align*}
These vector spaces carry the same structures as $(HH^*(X), HH_*(X))$,
namely $HT^*(X)$ is a ring, with multiplication induced by the
exterior product on polyvector fields; $H\Omega_*(X)$ is a module over
$HT^*(X)$, via contraction of polyvector fields with forms; and in
Section~\ref{sec:harMukprod} we shall define a pairing on
$H\Omega_*(X)$ which is a modification of the usual pairing of forms
given by cup product and integration on $X$.  (This modified inner
product is a more concrete generalization of the Mukai product
in~\cite{MukK3}.)

\subsection{}
In Section~\ref{sec:prel} we explain how to associate to an
integral transform $\Phi$ a map of graded vector spaces
\[ \Phi_*:H\Omega_*(X) \ra H\Omega_*(Y) \]
and we prove in Section~\ref{sec:harMukprod} that this association
satisfies the same adjointness properties as the similar association for
Hochschild homology discussed above.

\subsection{}
The connection between the Hochschild and harmonic structures is given
by the Hoch\-schild-Kostant-Rosenberg (HKR) isomorphism, which in
modern language can be written as a specific quasi-isomorphism
\[ I:\Delta^* \cO_\Delta \stackrel{\sim}{\lra} \bigoplus_i
\Omega_X^i[i], \] 
where $\Delta^*$ is the left derived functor of the usual pull-back
functor, and the right hand side of the quasi-isomorphism is the complex
which has $\Omega_X^i$ in the $-i$-th position, and all differentials
are zero.  The isomorphism $I$ induces isomorphisms of graded vector
spaces (Corollary~\ref{cor:ihkr})
\begin{align*}
I^\HKR & : HH^*(X) \stackrel{\sim}{\lra} HT^*(X), \\
I_\HKR & : HH_*(X) \stackrel{\sim}{\lra} H\Omega_*(X).
\end{align*}

\vspace{1mm}
\noindent
{\bf Theorem~\ref{thm:chernhar}.}  
{\em The composition
\[
\begin{diagram}[height=2em,width=2em,labelstyle=\scriptstyle]
K_0(X) & \rTo^{\ch} & HH_0(X) & \rTo^{I_\HKR} & \bigoplus_i H^i(X,
\Omega_X^i) 
\end{diagram}
\]
agrees with the usual Chern character map.}
\vspace{1mm}

This result was originally stated without proof and in an incomplete
form in a preprint by Markarian~\cite{Mar}.

As part of our proof of this theorem we prove the following result,
which provides an interesting interpretation of the Atiyah class in
view of the HKR isomorphism:

\vspace{2mm}
\noindent
{\bf Proposition~\ref{prop:etaatyiah}.}  
{\em The exponential of the universal Atiyah class is precisely the
map
\[ \cO_\Delta\stackrel{\eta}{\lra}\Delta_*\Delta^*\cO_\Delta
\stackrel{\Delta_*I}{\lra} \bigoplus_i \Delta_*\Omega_X^i[i], \]
where $\eta$ is the unit of the adjunction $\Delta^*\adjoint
\Delta_*$. }

\subsection{}
While the HKR isomorphism is well-behaved with respect to the Chern
character (in fact one can take Theorem~\ref{thm:chernhar} as a
definition of the differential forms-valued Chern character), it was
argued by Kontsevich~\cite{Kon} and Shoikhet~\cite{Sho} that $I^\HKR$,
$I_\HKR$ do not respect the Hochschild and harmonic structures.
Specifically, $I^\HKR$ is {\em not} a ring isomorphism.  However,
Kontsevich argued that as a consequence of his proof of the formality
conjecture, modifying $I^\HKR$ by the square root of the Todd genus
does in fact yield a ring isomorphism.  More precisely, denote by
$I^K$ the isomorphism
\[
\begin{diagram}[height=2em,width=2em,labelstyle=\scriptstyle]
I^K:HH^*(X) & \rTo^{I^\HKR} & HT^*(X) & \rTo^{\vee \Td_X^{-1/2}} &
HT^*(X).
\end{diagram}
\]
where the second map is given by the contraction of a polyvector field
with $\Td_X^{-1/2}$.  Then $I^K$ is a ring isomorphism~\cite[Claim
  8.4]{Kon}.

\subsection{}
A similar phenomenon can be seen on the level of homology theories:
the Mukai product that we define in~(\ref{subsec:mukprod}) does {\em
not} satisfy
\[ \langle \ch(\cE), \ch(\cF) \rangle = \chi(\cE, \cF) \]
as would have been expected from the similar property of Hochschild
homology.  The correct statement (already known to Mukai in the case
of K3 surfaces) is that
\[ \langle v(\cE), v(\cF) \rangle = \chi(\cE, \cF), \]
where 
\[ v(\cE) = \ch(\cE)\wedge\Td(X)^{1/2}. \]
These observations lead to the following conjecture:

\vspace{2mm}
\noindent
{\bf Conjecture~\ref{conj:mainconj}.}  
{\em
The maps
\[ I^K: HH^*(X) \ra HT^*(X), \quad\quad I_K: HH_*(X) \ra H\Omega_*(X),
\]
where $I^K$ is the composition
\begin{align*}
\begin{diagram}[height=2em,width=2em,labelstyle=\scriptstyle]
I^K:HH^*(X) & \rTo^{I^\HKR} & HT^*(X) & \rTo^{\vee \Td_X^{-1/2}} &
HT^*(X),
\end{diagram}
\intertext{and $I_K$ is given by} 
\begin{diagram}[height=2em,width=2em,labelstyle=\scriptstyle]
I_K:HH_*(X) & \rTo^{I_\HKR} & H\Omega_*(X) & \rTo^{\wedge \Td_X^{1/2}} & H\Omega_*(X),
\end{diagram}
\end{align*}
induce an isomorphism between the Hochschild and the harmonic
structures of $X$.  Concretely, $I^K$ is a ring isomorphism, $I_K$ is
an isometry with respect to the generalized Mukai product, and the two
isomorphisms are compatible with the module structures on
$H\Omega_*(X)$ and $HH_*(X)$, respectively.}
\vspace{2mm}

It is worthwhile observing that both $I^K$ and $I_K$ arise from the
same modification of the HKR isomorphism $I$~(\ref{subsec:samemod}).

\subsection{}
The main reason these results are interesting is because it has been
conjectured by Kontsevich~\cite{KonICM} that, in the case of a
Calabi-Yau manifold, $HH^*(X)$ should be the same as the ordinary
cohomology ring $H^*(X,\C)$ of the {\em mirror} $\check{X}$ of $X$.
In ~\cite{CalHH3} we shall expand this idea further, introducing a product
structure on the Hochschild homology of a Calabi-Yau orbifold and
arguing that its properties make it a good candidate for the mirror of
Chen-Ruan's~\cite{ChenRuan} orbifold cohomology theory.

Another application of the results in this paper, also to appear
in~\cite{CalHH3}, is a conceptual explanation of the results of the
computations of Fantechi and G\"ottsche~\cite{FanGot}, which show that
the orbifold cohomology of a symmetric product of abelian or K3
surfaces agrees with the cohomology of the Hilbert scheme of points on
the surface.  This explanation is a combination of the main result of
Bridgeland, King and Reid~\cite{BKR} with ideas of Verbitsky~\cite{Ver}
and with the derived category invariance of the Hochschild structure.

\subsection{}
The paper is structured as follows: after an introductory section in
which we discuss integral transforms and natural transformations
between them, we turn in Section~\ref{sec:harMukprod} to a definition
of the Mukai pairing on forms and to proofs of its basic functoriality
and adjointness properties.  Section~\ref{sec:hkr} is devoted to a
discussion of the HKR isomorphism and of the compatibility between the
Chern character~(\ref{subsec:chernchar}) and the usual one.  We
conclude with a discussion of the main conjecture and of possible ways
of proving it in Section~\ref{sec:mainconj}.

\vspace{1.5mm}
\noindent
\textbf{Acknowledgments.}  I have greatly benefited from conversations
with Tom Bridgeland, Andrew Kresch, Tony Pantev and Jonathan Block.  I
am thankful to Amnon Yekutieli for pointing out a mistake in an
earlier version of the paper.  The author's work has been supported by
an NSF postdoctoral fellowship and by travel grants and hospitality
from the University of Pennsylvania, the University of Salamanca,
Spain, and the Newton Institute in Cambridge, England.

\vspace{1.5mm}
\noindent
\textbf{Conventions.}  All the spaces involved are smooth algebraic
varieties proper over $\C$ (or any algebraic closed field of
characteristic zero), or compact complex manifolds.  We shall always
omit the symbols $\Ld$ and $\R$ in front of push-forward, pull-back
and tensor functors, but we shall consider them as derived except
where explicitly stated otherwise.  We shall write $\cF\otimes \mu$
where $\cF$ is a sheaf and $\mu$ is a morphism and mean by this the
morphism $1_\cF\otimes \mu$.  We shall use either $\wedge$ or $.$ for
the usual product in cohomology.  Serre duality notations and conventions
are presented in detail in Section~\ref{sec:prel}.

\section{Preliminaries}
\label{sec:prel}

In this section we provide a brief introduction to integral functors
on the level of derived categories and rational cohomology.  The
concepts and results are not new, dating back to Mukai's seminal
papers~\cite{MukAb},~\cite{MukK3}.  We also include several results
on traces and duality theory that will be needed later on.  

\subsection{}
\label{subsec:inttransf}
Let $X$ and $Y$ be complex manifolds, and let $\cE$ be an object in
$\D(X\times Y)$.  If $\pi_X$ and $\pi_Y$ are the projections from
$X\times Y$ to $X$ and $Y$, respectively, we define the integral
transform with kernel $\cE$ to be the functor
\[ \FMXY^\cE:\D(X) \ra \D(Y) \quad \FMXY^\cE(\scdot) = \pi_{Y, *}
(\pi_X^*(\scdot) \otimes \cE) \]
Likewise, if $\mu$ is any element of the ring $H^*(X\times Y, \Q)$,
we define the map
\[ \fmXY^\mu:H^*(X, \Q) \ra H^*(Y, \Q) \quad \fmXY^\mu(\scdot) =
\pi_{Y,*} (\pi_X^*(\scdot).\mu) \] 
and call it the integral transforms (in cohomology) associated to
$\mu$.

\subsection{}
\label{subsec:nattrans}
The association between objects of $\D(X\times Y)$ and integral
transforms is functorial: given a morphism $\mu:\cE\ra \cF$ between
objects of $\D(X\times Y)$, there is an obvious natural transformation
\[ \FMXY^\mu:\FMXY^\cE\Rightarrow \FMXY^\cF \]
given by
\[ \FMXY^\mu(\scdot) = \pi_{Y,*}(\pi_X^*(\scdot) \otimes \mu). \]

\subsection{}
There is a natural map between the derived category and the cohomology
ring, namely the exponential Chern character, $\ch:\D(X)\ra H^*(X,
\Q)$.  It commutes with pull-backs, and transforms tensor products
into cup products.  In an ideal world, it would also commute with
push-forwards, and then the diagram
\[
\begin{diagram}[height=2em,width=2em,labelstyle=\scriptstyle]
\D(X) & \rTo^{\FMXY^\cE} & \D(Y) \\
\dTo^{\ch} & & \dTo^{\ch} \\
H^*(X,\Q) & \rTo^{\fmXY^{\ch(\cE)}} & H^*(Y, \Q).
\end{diagram}
\]
would commute.  However, the Grothendieck-Riemann-Roch formula tells
us that we need to correct the commutation of push-forward and $\ch$
by the Todd classes of the spaces involved; more precisely, if
$\pi:X\ra Y$ is a locally complete intersection morphism, then
\[ \pi_*(\ch(\scdot).\Td(X)) = \ch(\pi_*(\scdot)).\Td(Y). \]

It is easy to see that there exists a unique formal series expansion
$\sqrt{1+c_1+c_2+\ldots}$ in the symbols $c_1, c_2, \ldots$, such that
\begin{align*}
\sqrt{1} & = 1\\
\sqrt{\mu.\nu} & = \sqrt{\mu}.\sqrt{\nu} \\
\intertext{and}
(\sqrt{\mu})^2 & = \mu
\end{align*}
for every space $X$ and any $\mu, \nu \in H^\even(X, \Q)$ with
constant term 1.  Its first three terms are
\[ \sqrt{1+c_1+c_2+\ldots} = 1 + \frac{1}{2} c_1 + \frac{1}{8}(4 c_2 -c_1^2) +
\frac{1}{16}(8c_3-4c_1 c_2+c_1^3) + \ldots\,. \]
Since the Todd class of any space $X$ is a sum of even cohomology
classes with constant term 1, $\sqrt{\Td(X)}$ is well defined by the
formula above, and we can define
\[ v:\D(X) \ra H^*(X, \Q) \quad v(\scdot) = \ch(\scdot).\sqrt{\Td(X)}.
\]
For an element $\cE$ of $\D(X)$ (on some space $X$), $v(\cE)$ will
be called the {\em Mukai vector} of $\cE$.  A straightforward
calculation shows that the diagram
\[
\begin{diagram}[height=2em,width=2em,labelstyle=\scriptstyle]
\D(X) & \rTo^{\FMXY^\cE} & \D(Y) \\
\dTo^{v} & & \dTo^{v} \\
H^*(X,\Q) & \rTo^{\fmXY^{v(\cE)}} & H^*(Y, \Q).
\end{diagram}
\]
commutes.  (This is a direct analogue of~\cite[Theorem~7.1]{CalHH1}.)
We shall denote the map $\fmXY^{v(\cE)}$ by $\Phi_*$, where $\Phi =
\FMXY^\cE$. 

\subsection{}
\label{subsec:composite}
Given complex manifolds $X, Y, Z$, and elements $\cE\in \D(X\times
Y)$ and $\cF\in \D(Y\times Z)$, define $\cF\circ\cE\in\D(X\times Z)$ by 
\[ \cF\circ \cE =  \pi_{XZ, *} (\pi_{XY}^* \cE \otimes
\pi_{YZ}^* \cF), \]
where $\pi_{XY}, \pi_{YZ}, \pi_{XZ}$ are the projections from $X\times
Y\times Z$ to $X\times Y$, $Y\times Z$ and $X\times Z$ respectively.
Similarly, if $\mu\in H^*(X\times Y, \Q)$, $\nu\in H^*(Y\times Z,
\Q)$, consider $\nu\circ\mu\in H^*(X\times Z, \Q)$ given by
\[ \nu\circ \mu = \pi_{XZ, *}(\pi_{XY}^* \mu . \pi_{YZ}^* \nu). \]
The reason behind the notation is the fact that
\begin{align*}
\FMYZ^\cF \circ \FMXY^\cE = \FMXZ^{\cF\circ\cE},
\intertext{and}
\fmYZ^\nu \circ \fmXY^\mu = \fmXZ^{\nu\circ \mu}.
\end{align*}
(The second result is standard; for a proof of the first one
see~\cite[1.4]{BonOrl}.)  Furthermore, it is a straightforward
calculation to check that
\[ v(\cF\circ\cE) = v(\cF) \circ v(\cE) \]
(\cite[3.1.10]{Cal}).  It follows that if $\Psi:\D(X)\ra \D(Y)$ and
$\Phi:\D(Y)\ra \D(Z)$ are integral transforms, then we have 
\[ (\Phi\circ \Psi)_* = \Phi_* \circ \Psi_* \]
(compare also to~\cite[Theorem~6.3]{CalHH1}).  Since it can be easily
checked that $\Id_* = \Id$, it follows that if $\Phi$ is an
equivalence of derived categories, then $\Phi_*$ is an isomorphism
$H^*(X, \Q) \ra H^*(Y, \Q)$.

\subsection{}
\label{subsec:vertHodge}
The map $\Phi_*$ does not respect the usual grading on the
cohomology rings of $X$ and $Y$, nor does it respect Hodge
decompositions.  However, it does respect the decomposition of
$H^*(X)$ by {\em columns} of the Hodge diamond: for every $i$,
$\Phi_*$ maps $H\Omega_i(X)$ to $H\Omega_i(Y)$,
\[ \Phi_* = \fmXY^{v(\cE)}: H\Omega_i(X) = \bigoplus_{q-p=i} H^{p,q}(X) 
\ra H\Omega_i(Y) = \bigoplus_{q-p=i} H^{p,q}(Y), \] 
because $v(\cE)$ consists only of classes of type $H^{p,p}(X\times
Y)$, and pushing-forward to $Y$ maps a class of type $(p,q)$ to a
class of type $(p-\dim Y, q-\dim Y)$.

This statement is the harmonic structure analogue of the fact that the
push-forward on Hochschild homology preserves the grading.

\section{The Mukai pairing on cohomology}
\label{sec:harMukprod}

In Section~\ref{sec:prel} we defined an isomorphism $\Phi_*: H^*(X,
\Q)\ra H^*(Y, \Q)$ associated to an equivalence of categories
$\Phi:\D(X)\ra \D(Y)$.  In the case of K3 surfaces, Mukai proved that
although $\Phi_*$ does not respect the usual intersection pairing on
the total cohomology rings of $X$ and of $Y$, it is an isometry with
respect to a modified version of this pairing.  He did this by showing
the more powerful result that maps on cohomology associated to adjoint
functors are themselves adjoint with respect to this modified pairing.
In this section we generalize this result to arbitrary complex
manifolds (not necessarily of dimension 2 or with trivial canonical
class), by defining a suitable generalization of Mukai's pairing.

\subsection{}
\label{subsec:elemprod}
The reason behind $\Phi_*$ being an isometry for the Mukai
product is the fact that an equivalence $\Phi:\D(X) \ra \D(Y)$ must
satisfy
\begin{align*}
\chi_X(\cF, \cG) & = \sum_{i} (-1)^i \dim \R \Hom^i(\cF, \cG) \\
	& = \sum_i (-1)^i \dim \R \Hom^i(\Phi\cF, \Phi\cG) \\
	& = \chi_Y(\Phi\cF, \Phi\cG).
\end{align*}
Thus, if we define a pairing on the algebraic part of $H^*(X, \Q)$ by
\[ \langle v(\cF), v(\cG)\rangle = \chi_X(\cF, \cG) \]
for all $\cE, \cF\in \D(X)$, then $\Phi_*$ is an isometry between the
algebraic subrings of $H^*(X, \Q)$ and $H^*(Y, \Q)$ (because $v$
commutes with $\Phi$).

\subsection{}
There are two problems with this definition: one is whether the above
pairing is well defined, another if we can extend it to a pairing on
the whole cohomology ring of $X$.  For K3 surfaces we have
\begin{align*}
\chi_X(\cF, \cG) & = \chi_X(\cF^\chk \otimes \cG) \\
	& = \int_X \ch(\cF^\chk).\ch(\cG).\Td(X) \\
	& = \int_X \ch(\cF^\chk).\sqrt{\Td(X)}.\ch(\cG).\sqrt{\Td(X)} \\
	& = \int_X v(\cF^\chk).v(\cG) \\
	& = \int_X v(\cF)^\chk.v(\cG),
\end{align*}
where $\cF^\chk = \R\Hom(\cF, \cO_X)$, and for a vector
\[ v=(v_0, v_2, v_4) \in H^0(X, \Q)\oplus H^2(X, \Q) \oplus H^4(X, \Q) \] 
$v^\chk$ is defined to equal $(v_0, -v_2, v_4)$.  Thus the pairing is
well defined in the K3 case (it only depends on the Mukai vectors of
$\cF$ and $\cG$, and not on $\cF$ and $\cG$ themselves).

\subsection{}
Our goal is to define $v^\chk$ for every $v\in H^\even(X, \Q)$ (and
eventually for every $v\in H^*(X, \Q)$), such that we have the
equality
\[ v(\cF^\chk) = v(\cF)^\chk, \]
which is the critical step in the previous computation.  When $X$ and
$Y$ are arbitrary complex manifolds the definition needs to take into
account the (possibly non-trivial) canonical class of $X$.

Define
\begin{align*}
& \tau : H^\even(X, \Q) \ra H^\even(X, \Q) 
\intertext{by}
& \tau(v_0, v_2, \ldots, v_{2n}) = (v_0, -v_2, v_4, \ldots, (-1)^n
v_{2n}).
\end{align*}
Then $\tau$ is easily checked to satisfy $\tau(v.w) =
\tau(v).\tau(w)$, and it is well known that $\ch(\cF^\chk) =
\tau(\ch(\cF))$.  Thus
\begin{align*}
v(\cF^\chk) & = \ch(\cF^\chk).\sqrt{\Td(X)} =
\tau(\ch(\cF)).\sqrt{\Td(X)} \\
	& = \tau\left(\frac{v(\cF)}{\sqrt{\Td(X)}}\right) \sqrt{\Td(X)} \\
	& = \tau(v(\cF)). \frac{1}{\sqrt{\ch(\omega_X)}},
\end{align*}
where the last equality is an immediate consequence of the formula
\[ \Td(T_X^\chk) = \Td(T_X).\exp(-c_1(T_X)) = \Td(T_X).\ch(\omega_X)
\]
(\cite[I.5.2]{FulLan}). 

Thus, if we define 
\[ v^\chk = \tau(v)\frac{1}{\sqrt{\ch(\omega_X)}} \]
for every $v\in H^\even(X, \Q)$, we have
\[ v(\cF^\chk) = v(\cF)^\chk \]
for all $\cF\in \D(X)$.

\subsection{}
To obtain a full generalization of the Mukai product we need to extend
the above mapping to all of $H^*(X, \Q)$.  A natural extension of the
involution $\tau$ is the map $\tau:H^*(X, \C) \ra H^*(X, \C)$ given by
\[ \tau(v_0, v_1, v_2, \ldots, v_{2n}) = (v_0, iv_1, -v_2, \ldots,
i^{2n}v_{2n}), \]
where $i = \sqrt{-1}$.  Its main properties are
\begin{enumerate}
\item $\tau(v.w) = \tau(v).\tau(w)$;
\item $\tau(\sqrt{v}) = \sqrt{\tau(v)}$ for any $v$ with leading term
equal to 1;
\item $\tau(\tau(v)) = v$ for any $v\in H^\even(X, \C)$;
\item $\tau(\ch(\cL)) = \ch(\cL^{-1}) = \ch(\cL)^{-1}$ for any line bundle 
$\cL$;
\item $\tau(f^*(v)) = f^*(\tau(v))$;
\item $f_*(\tau(v)) = (-1)^{\dim_\C X - \dim_\C Y} \tau (f_* v)$,
\end{enumerate}
where $f:X\ra Y$ is any proper morphism of complex manifolds.  The
proof of all these properties is immediate.

Thus, defining
\begin{align*}
\scdot^\chk & :H^*(X, \C) \ra H^*(X, \C)
\intertext{by}
v^\chk & = \tau(v).\frac{1}{\sqrt{\ch(\omega_X)}}
\end{align*}
extends in a natural way the operator $\scdot^\chk$ previously defined.

\subsection{}
\label{subsec:mukprod}
We can now tackle the generalized Mukai product:
\begin{definition}
\label{def:mukprodcoho}
Let $X$ be a complex manifold, and let $v,w \in H^*(X, \C)$.  Define
the product $\langle v,w\rangle$ by the formula
\[ \langle v, w\rangle = \int_X v^\chk.w, \]
where $v^\chk$ is defined above.  This product will be called the {\em
generalized Mukai product}.
\end{definition}

\subsection{}
It is interesting to compare this definition with a similar one that
appears in Hodge theory.  Define the Weyl operator, $\bar{\tau}$, by
$\bar{\tau}(v) = i^{p-q}v$ for $v\in H^{p,q}(X)$.  The pairing
\[ \langle v, w \rangle = \int_X \bar{\tau}(v). w \]
is the standard one that appears in the definition of a polarized
Hodge structure.  Observe that the analogy between the Mukai pairing
as a mirror to the usual Poincar\'e pairing holds, if we take this in
the sense of matching polarizations: the map $\tau$ is formally the
mirror of $\bar{\tau}$ (if we mirror the Hodge diamond, $\tau$ gets
transformed into $\bar{\tau}$).  We do not have a good understanding
of the $1/\sqrt{\ch(\omega_X)}$ term that appears in the definition of
$v^\chk$.

\begin{proposition}
\label{prop:fmadj}
Let $X$ and $Y$ be complex manifolds, and $\Phi:\D(Y)\ra\D(X)$ and
$\Psi:\D(X)\ra \D(Y)$ be adjoint integral transforms
($\Psi$ is a left adjoint to $\Phi$).  Then $\Psi_*$ is a left
adjoint to $\Phi_*$ with respect to the generalized Mukai
product; in other words, we have
\[ \langle v, \Phi_* w)\rangle_Y = \langle\Psi_*v, 
w\rangle_X \]
for all $v\in H^*(Y, \C)$, $w\in H^*(X, \C)$.
\end{proposition}

\begin{remark}
When $v$ and $w$ are Mukai vectors of elements in $\D(Y)$ and $\D(X)$,
the result is a trivial consequence of the discussion
in~(\ref{subsec:elemprod}).  The actual content is that the result holds
for all $v,w$.  
\end{remark}

\begin{corollary}
\label{cor:isometry}
Under the hypotheses of Proposition~\ref{prop:fmadj}, assume
furthermore that $\Phi$ is an equivalence of categories.  Then
$\Phi_*:H^*(X, \C) \ra H^*(Y, \C)$ is an isometry with
respect to the generalized Mukai product.
\end{corollary}

\begin{proof}
See the proof of~\cite[Corollary 7.5]{CalHH1}.
\end{proof}

\begin{proof}[Proof (of Proposition~\ref{prop:fmadj})]
Assume $\Phi = \FMXY^\cE$, and let $\cE^* = \cE^\chk \otimes
\pi_Y^*\omega_Y[\dim Y]$, so that $\Psi=\FMYX^{\cE^*}$.  Define $e =
v(\cE)$ and $e^* = v(\cE^*)$.  We have
\begin{align*}
e^* & = \tau\left(\frac{e}{\sqrt{\Td(X\times Y)}}\right).
\pi_Y^*(\ch(\omega_Y)).(-1)^{\dim Y}.\sqrt{\Td(X\times Y)} = \\
	& = (-1)^{\dim Y} \tau(e).\frac{\pi_Y^*\sqrt{\ch(\omega_Y)}}
		{\pi_X^*\sqrt{\ch(\omega_X)}}
\intertext{and thus}
\tau(e^*) & = (-1)^{\dim Y} e.
	\frac{\pi_X^*\sqrt{\ch(\omega_X)}}{\pi_Y^*\sqrt{\ch(\omega_Y)}}.
\end{align*}
We then have
\begin{align*}
\langle \Psi_* v, w \rangle & = \langle\fmYX^{e^*}(v),w\rangle = \int_X \fmYX^{e^*}(v)^\chk.w = \int_X
\tau(\fmYX^{e^*}(v)).\frac{1}{\sqrt{\ch(\omega_X)}}.w \\
	 & = \int_X
	 \tau(\pi_{X,*}(\pi_Y^*v.e^*)).\frac{1}{\sqrt{\ch(\omega_X)}}.w
	 \\
	& = (-1)^{\dim Y} \int_X \pi_{X,*}(\tau(\pi_Y^*
	 v).\tau(e^*)).\frac{1}{\sqrt{\ch(\omega_X)}}.w \\
	& = (-1)^{\dim Y} \int_{X\times Y} \tau(\pi_Y^*
	 v).\tau(e^*).\frac{1}{\pi_X^*\sqrt{\ch(\omega_X)}}.\pi_X^*w\\
	& = (-1)^{\dim Y} \int_{X\times Y} \tau(\pi_Y^* v).(-1)^{\dim Y}. e.
	\frac{\pi_X^*\sqrt{\ch(\omega_X)}}{\pi_Y^*\sqrt{\ch(\omega_Y)}}.\frac{1}{\pi_X^*\sqrt{\ch(\omega_X)}}.\pi_X^*w\\
	& = \int_{X\times Y}
	 \pi_Y^*(\tau(v)).\frac{1}{\pi_Y^*\sqrt{\ch(\omega_Y)}}.e.\pi_X^*w\\
	& = \int_Y
	 \tau(v).\frac{1}{\sqrt{\ch(\omega_Y)}}.\pi_{Y,*}(e.\pi_X^* w)
	 \\
	& = \int_Y v^\chk.\fmXY^e(w) = \langle v, \fmXY^e(w)\rangle \\
	& = \langle v, \Phi_* w \rangle.
\end{align*}
\end{proof}

\section{The Hochschild-Kostant-Rosenberg theorem and the Chern character}
\label{sec:hkr}

In this section we study the relationship between the Hochschild and
harmonic structures.  We provide a discussion of the connection
between the usual Chern character and the one introduced
in~\cite{CalHH1}.

\subsection{}
The starting point of our analysis is the following theorem:
\begin{theorem}[Hochschild-Kostant-Rosenberg~\cite{HKR},
    Kontsevich~\cite{Kon}, Swan~\cite{Swa}, Yekutieli~\cite{Yek}]
\label{thm:hkr}
Let $X$ be a smooth, quasi-projec\-tive variety of dimension $n$,
and let $\Delta:X\ra X\times X$ be the diagonal embedding.  Then there
exists a quasi-isomorphism
\[ I:\Delta^*\cO_\Delta \stackrel{\sim}{\lra}
\bigoplus_i \Omega_X^i[i], \]
where the right hand side denotes the complex whose $-i$-th term is
$\Omega_X^i$, and all differentials are zero.
\end{theorem}

\begin{proof}
(This is nothing but a brief recounting of the results in~\cite{Yek},
and the reader should consult [loc.cit.] for more details.)  Recall
that if $R$ is a commutative $\C$-algebra there exists a standard
resolution of $R$ as an $R^e = R\otimes_\C R$-module.  For $i\geq 0$
let
\[ \cB_i(R) = R^{\otimes (i+2)}, \]
where the tensor product is taken over $\C$.  It is an $R^e$-module by
multiplication in the first and last factor.  The bar resolution is
defined to be the complex of $R^e$-modules
\[ \cdots\ra \cB_i(R) \ra \cdots \ra \cB_1(R) \ra \cB_0(R) \ra 0, \]
with differential
\begin{eqnarray*} 
\lefteqn{d(a_0\otimes a_1\otimes\cdots\otimes a_i) = } \\
& & a_0a_1 \otimes a_2 \otimes\cdots \otimes a_i \,-\, 
a_0\otimes a_1a_2 \otimes \cdots \otimes a_i \,+\, \cdots \,+\, \\ 
& & (-1)^{i-1}a_0\otimes a_1\otimes\cdots \otimes
a_{i-1}a_i. 
\end{eqnarray*}
It is an exact complex, except at the last step where the cohomology
is $R$. Thus it is a resolution of $R$ in
$R^e$-$\gMod$~\cite[1.1.12]{Lod}.

If $X$ were affine, $X=\Spec R$, we could use the above resolution to
compute $\Delta^* \cO_\Delta$: indeed, $\cO_\Delta$ is nothing but $R$
viewed as an $R^e = \cO_{X\times X}$-module, and the modules
$\cB_i$ are $R^e$-flat.  The complex obtained by tensoring the
bar resolution over $R^e$ with $R$ is called the bar complex:
\[ \cdots\ra \cC_i(R) \ra \cdots \ra \cC_1(R) \ra \cC_0(R) \ra 0, \]
where 
\[ \cC_i(R) = \cB_i(R) \otimes_{R^e} R, \]
and the differential is obtained from the differential of $\cB_\cdot(R)$.

Problems arise when one tries to sheafify the bar resolution to
obtain a complex of sheaves on a scheme: the resulting sheaves are
ill-behaved (in particular, not quasi-coherent).
As a replacement, Yekutieli proposed to used the complete bar
resolution, which he defined in~\cite{Yek}.  For $i\geq 0$, let
$\sX^i$ be the formal completion of the scheme $X^i = X\times
\cdots\times X$ along the small diagonal.  Define
\[ \hcB_i(X) = \cO_{\sX^{i+2}}, \]
which is a sheaf of abelian groups on the topological space $X$.
Yekutieli argued that one can formally complete and sheafify the
original bar resolution to get the complete bar resolution
\[ \cdots\ra \hcB_i(X) \ra \cdots \ra \hcB_1(X) \ra \hcB_0(X) \ra 0, \]
where the maps are locally obtained from the maps of the original bar
complex, by noting that these are continuous for the topologies with
respect to which we are completing.  The complete bar resolution is an
exact resolution of $\cO_\Delta$ by sheaves of flat $\cO_{X\times
X}$-modules (see remark following Proposition 1.4 and proof of
Proposition 1.5 in~\cite{Yek}).  Over an affine open set $U = \Spec R$
of $X$, $\Gamma(U, \hcB_i(X))$ is the completion $\hcB_i(R)$ of
$\cB_i(R)$ at the ideal $I_i$ which is the kernel of the
multiplication map $\cB_i(R) = R^{\otimes i} \ra R$.

One can take the complete bar resolution as a flat resolution of
$\cO_\Delta$ on $X\times X$, and use it to compute
$\Delta^*\cO_\Delta$.  This is the same as tensoring the complete bar
resolution over $\cO_{X\times X}$ with $\cO_\Delta$.  The resulting
complex is called the complex of complete Hochschild chains of $X$
(see~\cite[Definition 1.3]{Yek} for details),
\[ \cdots\ra \hcC_i(X) \ra \cdots \ra \hcC_1(X) \ra \hcC_0(X) \ra 0, \]
where 
\[ \hcC_i(X) = \hcB_i(X) \otimes_{\cO_{X\times X}} \cO_\Delta. \]
Over an affine open set $U=\Spec R$, $\Gamma(U, \hcC_i(X))$ is the
completion $\hcC_i(R)$ of $\cC_i(R)$ at $I_i$ (as a $\cB_i(R)$-module).

Over any affine open $U= \Spec R$ define
\[ I_i:\cC_i(R) \ra \Omega_{R/k}^i \]
by setting
\[ I_i((1\otimes a_1\otimes\cdots\otimes a_i \otimes 1) \otimes_{R^e}
1) = d a_1\wedge da_2 \wedge\cdots \wedge da_i. \] 
These maps are continuous with respect to the topology that is used
for completing~(\cite[Lemma 4.1]{Yek}), so they can be completed and
sheafified to maps
\[ I_i:\hcC_i(X) \ra \Omega_X^i. \]
They also commute with the zero differentials of the complex $\oplus_i
\Omega_X^i$, so they assemble to a morphism of complexes
\[ I:\Delta^*\cO_\Delta \ra \bigoplus_i \Omega_X^i[i] \]
which can be seen to be a quasi-isomorphism in characteristic
0~(\cite[Theorem 4.6.1.1]{Kon},~\cite[Proposition 4.4]{Yek}).  In the
affine case this is essentially the Hochschild-Kostant-Rosenberg
theorem~\cite{HKR}.
\end{proof}

\begin{corollary}
\label{cor:ihkr}
The Hochschild-Kostant-Rosenberg isomorphism $I$ induces isomorphisms
of graded vector spaces
\begin{align*}
I^\HKR & : HH^*(X) \stackrel{\sim}{\lra} HT^*(X), \\
I_\HKR & : HH_*(X) \stackrel{\sim}{\lra} H\Omega_*(X).
\end{align*}
\end{corollary}

\begin{proof}
\begin{align*}
HH^k(X) & = \Hom_{X\times X}(\cO_\Delta, \cO_\Delta[k]) \iso
\Hom_X(\Delta^*\cO_\Delta, \cO_X[k]) \\
& \iso \Hom_X(\bigoplus_i \Omega_X^i[i], \cO_X[k]) = \bigoplus_i
H^{k-i}(X, \bigwedge^i T_X) = HT^k(X),
\intertext{and}
HH_k(X) & = \Hom_{X\times X}(\Delta_!\cO_X[k],\cO_\Delta) \iso
\Hom_X(\cO_X[k], \Delta^*\cO_\Delta) \\
& \iso \Hom_X(\cO_X[k], \bigoplus_i \Omega_X^i[i]) = \bigoplus_i
H^{i-k}(X, \Omega_X^i) = H\Omega_k(X).
\end{align*}
\end{proof}

\subsection{}
We are now interested in understanding how the above isomorphisms
relate the Chern character $K_0(X) \ra HH_0(X)$ defined in the
introduction to the usual Chern character.

Let $\Omega_\Delta^{\otimes i}$ and $\Omega_\Delta^i$ denote the
push-forwards by $\Delta$ of $\Omega_X^{\otimes i}$ and $\Omega_X^i$,
respectively.  (Here the tensor product is taken over $\cO_X$.)  Let 
\[ \epsilon:\Omega_X^{\otimes i} \ra \Omega_X^i \]
be the antisymmetrization map which acts as
\[ v_1\otimes v_2\otimes\cdots v_i \mapsto \frac{1}{i!}
\sum_{\sigma\in\Sigma_i} (-1)^{\epsilon(\sigma)} v_{\sigma_1}\otimes
v_{\sigma_2}\otimes\cdots v_{\sigma_i}. \]
By an abuse of notation, we shall also denote by $\epsilon$ the
push-forward
\[ \Delta_*\epsilon:\Omega_\Delta^{\otimes i} \ra \Omega_\Delta^i. \]

\begin{definition}
Define the universal Atiyah class to be the class
\[ \alpha_1\in\Ext^1_{X\times X}(\cO_\Delta, \Omega_\Delta^1), \]
of the extension
\[ 0\ra \Omega_\Delta^1 \ra \cO_{\Delta^{(2)}} \ra \cO_\Delta \ra 0,
\]
where $\cO_{\Delta^{(2)}}$ is the second infinitesimal neighborhood of
the diagonal in $X\times X$.  Furthermore, define $\alpha_i$ for
$i\geq 0$ by the formula
\[ \alpha_i = \epsilon \circ (\pi_2^*\Omega_X^{\otimes (i-1)}
\otimes \alpha_1) \circ (\pi_2^*\Omega_X^{\otimes (i-2)} \otimes
\alpha_1) \circ \cdots\circ \alpha_1:\cO_\Delta\ra \Omega_\Delta^i[i]. \] 
The exponential Atiyah class $\exp(\alpha)$ is defined by the formula
\[ \exp(\alpha) =
1+\alpha_1+\alpha_2+\cdots+\alpha_n:\cO_\Delta\ra\bigoplus_i
\Delta_* \Omega_X^i[i]. \]
\end{definition}

This definition requires a short explanation.  Recall that given an
object $\cE\in\D(X)$, the Atiyah class of $\cE$ is the class
\[ \alpha_1(\cE)\in\Ext^1_X(\cE, \cE\otimes \Omega_X^1) \]
of the extension on $X$
\[ 0\ra \cE\otimes \Omega_X^1\ra J^1(\cE)\ra \cE \ra 0 \]
where $J^1(\cE)$ is the first jet bundle of $\cE$~\cite[1.1]{KapAt}.
A natural way to construct this extension is to consider the natural
transformation $\FMXX^{\alpha_1}$ associated to the universal Atiyah
class
\[ \alpha_1:\cO_\Delta\ra \Omega_\Delta^1[1] \]
between the identity functor and the ``tensor by $\Omega_X^1[1]$''
functor.  The value $\FMXX^{\alpha_1}(\cE)$ of this natural
transformation on $\cE$ is precisely the Atiyah class $\alpha_1(\cE)$
of $\cE$ (see, for example,~\cite[10.1.5]{BlueBook}).  The $i$-th
component of the Chern character of $\cE$ is then obtained as
\[ \ch_i(\cE) = \Tr_\cE(\alpha_i(\cE)) \]
where 
\[ \alpha_i(\cE) = \epsilon\circ
(\Omega_X^{\otimes(i-1)}\otimes\alpha_1(\cE)) \circ
(\Omega_X^{\otimes(i-2)}\otimes\alpha_1(\cE)) \circ \cdots \circ
\alpha_1(\cE):\cE\ra \cE\otimes \Omega_X^i[i]. \]
(See~\cite[10.1.6]{BlueBook} for details.)  Our definition of
$\alpha_i:\cO_\Delta\ra \Delta_*\Omega_X^i[i]$ has been tailored to
mimic this definition: $\alpha_i(\cE)$ will be precisely the value on
$\cE$ of the natural transformation associated to the morphism
$\alpha_i$.  Therefore, if we consider the natural transformation
$\FMXX^{\exp(\alpha)}$ associated to $\exp(\alpha)$, its value
\[ \Phi^{\exp(\alpha)}(\cE):\cE\ra\bigoplus_i \cE\otimes
\Omega_X^i[i] \]
on $\cE$ will satisfy
\[ \ch_{\mathrm{orig}}(\cE) = \Tr_\cE(\Phi^{\exp(\alpha)}(\cE)), \]
where $\ch_{\mathrm{orig}}(\cE)$ is the usual Chern character of
$\cE$.

\begin{proposition}
\label{prop:etaatyiah}
The exponential $\exp(\alpha)$ of the universal Atiyah class is
precisely the map
\[ \cO_\Delta\stackrel{\eta}{\lra}\Delta_*\Delta^*\cO_\Delta
\stackrel{\Delta_*I}{\lra} \bigoplus_i \Delta_*\Omega_X^i[i], \]
where $\eta$ is the unit of the adjunction $\Delta^*\adjoint \Delta_*$.
\end{proposition}

\begin{proof}
We divide the proof of this proposition into several steps, to make it
more manageable.  We will use the notations used in the proof of
Theorem~\ref{thm:hkr}. 
\vspace{2mm}

\noindent {\em Step 1.} 
Consider the exact sequence
\[ 0\ra \Omega_\Delta^1 \ra \cO_{\Delta^{(2)}} \ra \cO_\Delta \ra 0 \]
which defines the universal Atiyah class $\alpha_1$.  Tensoring it by
the locally free sheaf $\pi_2^* \Omega_X^{\otimes i}$ yields the exact
sequence
\[ 0 \ra \Omega_\Delta^{\otimes (i+1)} \ra \cO_{\Delta^{(2)}} \otimes 
\pi_2^* \Omega_X^{\otimes i}\ra \Omega_\Delta^{\otimes i} \ra 0. \]
Stringing together these exact sequences for successive values of $i$ we
construct the exact sequence
\[ 0\ra \Omega_\Delta^{\otimes i} \ra \cO_{\Delta^{(2)}} \otimes
\pi_2^*\Omega_X^{\otimes (i-1)} \ra \cO_{\Delta^{(2)}} \otimes
\pi_2^*\Omega_X^{\otimes (i-2)} \ra \cdots \ra \cO_{\Delta^{(2)}}\ra
\cO_\Delta \ra 0, \]
whose extension class is precisely
\[ (\pi_2^*\Omega_X^{\otimes (i-1)}
\otimes \alpha_1) \circ (\pi_2^*\Omega_X^{\otimes (i-2)} \otimes
\alpha_1) \circ \cdots\circ \alpha_1:\cO_\Delta\ra
\Omega_\Delta^{\otimes i}[i]. \]

\noindent {\em Step 2.}
We claim that there exists a map $\phi_\cdot$ of exact sequences
\begin{diagram}[height=2em,width=2em,labelstyle=\scriptstyle]
\cdots & \rTo & \hcB_i(X) & \rTo & \hcB_{i-1}(X) &
\rTo & \cdots & \rTo & \hcB_0(X) & \rTo & \cO_\Delta & \rTo & 0 \\
& & \dTo_{\phi_i'} & & \dTo_{\phi_{i-1}} & & & & \dTo_{\phi_0} & & \dEqual
& & \\
0 & \rTo & \Omega_\Delta^{\otimes i} & \rTo & \cO_{\Delta^{(2)}}
\otimes \pi_2^*\Omega_X^{\otimes (i-1)} & \rTo & \cdots
& \rTo & \cO_{\Delta^{(2)}} & \rTo & \cO_\Delta & \rTo & 0,
\end{diagram}
where the top row is the (augmented) completed bar resolution defined
in the proof of Theorem~\ref{thm:hkr}, and the bottom row is the one
defined in Step 1.  It is sufficient to define the maps in a local
patch $U=\Spec R$.  Let $I = I_2 = \ker(R\otimes R \ra R)$ be the
ideal defining the diagonal in $U\times U$.  Consider the maps
\[ \phi_i:\cB_i(R) = R^{\otimes (i+2)} \ra (R\otimes R)/I^2 \otimes_R
\Omega_R^{\otimes_R i} \]  
defined by 
\[ \phi_i(a_0\otimes a_1\otimes \cdots \otimes a_{i+1}) = (a_0\otimes
a_{i+1} + I^2) \otimes_R da_1 \otimes_R da_2
\otimes_R \cdots \otimes_R da_i \] 
(we write $\Omega_R$ on the right because we use $\pi_2^*$).  The
same argument as the one in the proof of~\cite[Lemma 4.1]{Yek} shows
that these maps are continuous with respect to the adic topology used
to complete $\cB_i(R) = R^{\otimes (i+2)}$, thus the maps $\phi_i$
descend to maps
\[ \phi_i:\hcB_i(R) \ra (R\otimes R)/I^2 \otimes_R
\Omega_R^{\otimes_R i} \]
which sheafify to give the desired maps 
\[ \phi_i: \hcB_i(X) \ra \cO_{\Delta^{(2)}} \otimes \pi_2^*
\Omega_X^{\otimes i}. \]
The map $\phi_i'$ is the composition 
\[ \hcB_i(X) \stackrel{\phi_i}{\lra} \cO_{\Delta^{(2)}}
\otimes \pi_2^* \Omega_X^{\otimes i} \ra \cO_\Delta \otimes \pi_2^*
\Omega_X^{\otimes i} = \Omega_\Delta^{\otimes i}. \]
\vspace{2mm}

\noindent {\em Step 3.}  
We now need to check the commutativity of the squares in the above
diagram.  Note that since everything is local, we can assume we are in
an open patch $U=\Spec R$, $U\times U = \Spec R\otimes R$.  The ideal
$I$ in $R\otimes R$ is generated by expressions of the form $r\otimes 1
-1\otimes r$ for $r\in R$.  Then a relevant square in the above
diagram (before completing) is
\[
\begin{diagram}[height=2em,width=2em,labelstyle=\scriptstyle]
R\otimes R\otimes R\otimes R & \rTo^{d_1} & R\otimes R\otimes R \\
\dTo_{\phi_2} & & \dTo_{\phi_1} \\ 
(R\otimes R)/I^2 \otimes_R I/I^2 \otimes_R I/I^2 & \rTo^{d_1'} &
(R\otimes R)/I^2 \otimes_R I/I^2, 
\end{diagram}
\]
where $(R\otimes R)/I^2$ is considered a right $R$-module by
multiplication in the second factor, and $I/I^2$ is considered an
$R$-module by multiplication in either factor (the two module
structures are the same).  The maps in this diagram are:
\begin{align*}
& d_1(1\otimes b\otimes c\otimes 1) = b\otimes c\otimes 1 - 1\otimes
bc\otimes 1 + 1\otimes b\otimes c,\quad\mbox{the Hochschild
differential} \\ 
& d_1'((r+I^2)\otimes_R (i+I^2)\otimes_R (i'+I^2)) = (ri+I^2)\otimes_R
(i'+I^2), \\
& \phi_1(a\otimes b\otimes c) = (a\otimes c + I^2)\otimes_R (b\otimes
1 - 1\otimes b + I^2), \\
& \phi_2(1\otimes b\otimes c\otimes 1) = (1\otimes 1 + I^2)\otimes_R
(b\otimes 1 - 1\otimes b + I^2) \otimes_R (c\otimes 1 - 1\otimes c +
I^2).
\end{align*}
Omitting the $+I^2$ terms, we have
\begin{align*}
\phi_1(d_1 & (1\otimes b\otimes c\otimes 1)) =
\phi_1(b\otimes c\otimes
1 - 1\otimes bc\otimes 1 + 1\otimes b\otimes c) \\
& = (b\otimes 1)\otimes_R (c\otimes 1 - 1\otimes c) - (1\otimes
1)\otimes_R (bc\otimes 1-1\otimes bc) + (1\otimes c)\otimes_R (b\otimes
1 - 1\otimes b) \\
\intertext{which, using the right module structure on $R\otimes R$, equals}
& = (b\otimes 1)\otimes_R (c\otimes 1 - 1\otimes c) - (1\otimes
1)\otimes_R (bc\otimes 1-1\otimes bc) + (1\otimes
1)\otimes_R (b\otimes c - 1\otimes bc) \\
& = (b\otimes 1)\otimes_R (c\otimes 1 - 1\otimes c) - (1\otimes
1)\otimes_R (bc\otimes 1 - b\otimes c) \\
\intertext{and note that, since the last term is equal to
  $1\otimes bc - c\otimes b$ modulo $I^2$,}
& = (b\otimes 1)\otimes_R (c\otimes 1 - 1\otimes c) - (1\otimes
1)\otimes_R (1\otimes bc - c\otimes b) \\
& = (b\otimes 1)\otimes_R (c\otimes 1 - 1\otimes c) - (1\otimes
b)\otimes_R (1\otimes c - c\otimes 1) \\
& = (b\otimes 1 - 1\otimes b) \otimes_R (c\otimes 1 - 1\otimes
c) \\
& = d_1'((1\otimes 1)\otimes_R (b\otimes 1 - 1\otimes b)
\otimes_R (c\otimes 1 - 1\otimes c)) \\
& = d_1'(\phi_2(1\otimes b\otimes c\otimes 1)). 
\end{align*}
Similar computations ensure the commutativity of the other squares.
\vspace{2mm}

\noindent {\em Step 4.}  
Observe that there exists a natural map $\eta$ from the bar resolution
$\hcB_\cdot(X)$ to the bar complex $\hcC_\cdot(X) = \hcB_\cdot(X)
\otimes_{X\times X} \cO_\Delta$, simply given by $1\otimes \mu$ where
$\mu:\cO_{X\times X} \ra \cO_\Delta$ is the natural projection.  This
map is immediately seen to be precisely the unit $\eta$ of the
adjunction $\Delta^*\adjoint \Delta_*$.

It is now obvious that the composite
\[ \hcB_i(X) \stackrel{\phi_i'}{\lra}
\Omega_\Delta^{\otimes i} \stackrel{\epsilon}{\lra} \Omega_\Delta^i \]
is precisely the same as the map
\[ \hcB_i(X) \stackrel{\eta_i}{\lra} \hcC_i(X)
\stackrel{\Delta_*I_i}{\lra} \Omega_\Delta^i, \] 
where $\eta_i$ is the $i$-th component of $\eta$, locally (before
completion) given by 
\[ a_0 \otimes a_1 \otimes \cdots \otimes a_{i+1} \mapsto a_0 a_{i+1}
\otimes a_1 \otimes \cdots \otimes a_i, \] 
and $\Delta_*I_i$ is the $i$-th component of the HKR isomorphism. 

Now chopping off at the last step the two exact sequences we have
studied above we get the diagram
\[ 
\begin{diagram}[height=2em,width=2em,labelstyle=\scriptstyle]
\cdots & \rTo & \hcB_i(X) &
\rTo & \hcB_{i-1}(X) & \rTo & \cdots & \rTo & \hcB_0(X) & \rTo & 0\\
 & & \dTo_{\phi_i'} & & \dTo_{\phi_{i-1}} & & & & 
\dTo_{\phi_0} & & \\
0 & \rTo & \Omega_\Delta^{\otimes i}
& \rTo & \cO_{\Delta^{(2)}} \otimes \pi_2^*\Omega_X^{\otimes (i-1)} & \rTo & \cdots
& \rTo & \cO_{\Delta^{(2)}} & \rTo & 0 \\
 & & \dTo_{p_i} & & & & & & & & & &\\
 & & \Omega_\Delta^{\otimes i} & & & & & & & & & & \\
 & & \dTo_{\epsilon} & & & & & & & & & & \\
 & & \Omega_\Delta^i, & & & & & & & & & &
\end{diagram}
\]
which can be thought of as a map from the top complex (which
represents $\cO_\Delta$) to $\Omega_\Delta^i[i]$.  In fact what we
have is a factoring
\[ \cO_\Delta \stackrel{p_i\circ\phi_\cdot}{\lra} \Omega_\Delta^{\otimes i}
\stackrel{\epsilon}{\lra} \Omega_\Delta^i \]
of the map 
\[ \epsilon \circ p_i \circ \phi_\cdot = \Delta_* I_i \circ \eta, \]
where $\phi_\cdot$ is the map of complexes appearing at the top of the
above diagram. However, note that both the source and the target of
$\phi_\cdot$ are naturally isomorphic (in $\D(X\times X)$) to
$\cO_\Delta$, and then $\phi_\cdot$ can be viewed as the identity map
$\cO_\Delta \ra \cO_\Delta$.  Under these identifications we conclude 
\[ \epsilon \circ p_i = \Delta_* I_i \circ \eta. \]

But the construction of $p_i$ is such that it is represented by the
$i$-step extension
\[
\begin{diagram}[height=2em,width=2em,labelstyle=\scriptstyle]
0 & \rTo & \Omega_\Delta^{\otimes i}
& \rTo & \cO_{\Delta^{(2)}} \otimes \pi_2^*\Omega_X^{\otimes (i-1)} & \rTo & \cdots
& \rTo & \cO_{\Delta^{(2)}} & \rTo &\cO_\Delta & \rTo & 0,
\end{diagram}
\]
whose class we argued is
\[ (\pi_2^*\Omega_X^{\otimes (i-1)}
\otimes \alpha_1) \circ (\pi_2^*\Omega_X^{\otimes (i-2)} \otimes
\alpha_1) \circ \cdots\circ \alpha_1:\cO_\Delta\ra
\Omega_\Delta^{\otimes i}[i]. \]
Therefore 
\[ p_i = (\pi_2^*\Omega_X^{\otimes (i-1)}
\otimes \alpha_1) \circ (\pi_2^*\Omega_X^{\otimes (i-2)} \otimes
\alpha_1) \circ \cdots\circ \alpha_1:\cO_\Delta\ra
\Omega_\Delta^{\otimes i}[i], \]
and hence
\[ \alpha_i = \epsilon \circ p_i  = \Delta_* I_i \circ \eta. \]
We conclude that
\[ \exp(\alpha) = \bigoplus_i \alpha_i  = \bigoplus \Delta_* I_i \circ
\eta = \Delta_*I \circ \eta. \]
\end{proof}

\begin{theorem}
\label{thm:chernhar}
The composition
\[
\begin{diagram}[height=2em,width=2em,labelstyle=\scriptstyle]
K_0(X) & \rTo^{\ch}& HH_0(X)&\rTo^{I_\HKR} \bigoplus_i 
H^i(X, \Omega^i_X)
\end{diagram} \]
is the usual Chern character map.
\end{theorem}

\begin{proof}
Let $\cF\in K_0(X)$, and let 
\[ \ch(\cF)\in HH_0(X) = \Hom_{X\times X}(\Delta_!\cO_{\Delta}, \cO_\Delta)
\]
be the Chern character defined in~(\ref{subsec:chernchar}).  Let 
\[ \ch'(\cF)\in \Hom_X(\cO_X, \Delta^*\cO_\Delta) \]
be the element that corresponds to $\ch(\cF)$ under the adjunction
$\Delta_!\adjoint \Delta^*$.  If $\mu'$ is any element of
$\Hom_X(\Delta^*\cO_\Delta, S_X)$ and 
\[ \mu = \Delta_*\mu'\circ \eta \]
is the corresponding element of $\Hom_{X\times X}(\cO_\Delta,
S_\Delta)$ under the adjunction $\Delta^*\adjoint\Delta_*$, the
construction of $\Delta_!$ is such that
\[ \Tr_X(\mu'\circ \ch'(\cF)) = \Tr_{X\times X}(\mu\circ \ch(\cF)). \]
(Here $\eta:\cO_\Delta \ra \Delta_*\Delta^*\cO_\Delta$ is the unit of
the adjunction.)

On the other hand, the definition of $\ch(\cF)$ is such that for any
$\mu$, 
\[ \Tr_{X\times X}(\mu\circ \ch(\cF)) =
\Tr_X(\pi_{2,*}(\pi_1^*\cF\otimes \mu)), \]
and $\ch(\cF)$ is the unique element in $HH_0(X)$ with this property.
We then have
\begin{align*}
\Tr_X(\mu'\circ \ch'(\cF)) & = \Tr_{X\times X}(\mu \circ \ch(\cF)) =
\Tr_X(\pi_{2,*}(\pi_1^*\cF\otimes \mu)) \\
& = \Tr_X(\pi_{2,*}(\pi_1^*\cF\otimes(\Delta_*\mu' \circ \eta))) \\
& = \Tr_X(\pi_{2,*}(\pi_1^*\cF\otimes \Delta_*\mu') \circ
\pi_{2,*}(\pi_1^*\cF\otimes \eta)) \\
& = \Tr_X(\cF\otimes \mu' \circ \Phi^\eta(\cF))\\
& = \Tr_X(\mu' \circ \Tr_\cF(\Phi^{\eta}(\cF))),
\end{align*}
where the last equality is~\cite[Lemma~2.4]{CalHH1}.
Since the trace induces a non-degenerate pairing and the above
equalities hold for any $\mu'$, it follows that
\[ \ch'(\cF) = \Tr_\cF(\Phi^{\eta}(\cF)). \]
Applying the isomorphism $I$ to both sides we conclude that 
\[ I_\HKR(\ch(\cF)) = I\circ\ch'(\cF)) = I\circ\Tr_\cF(\Phi^{\eta}(\cF)) =
\Tr_\cF(\Phi^{\exp(\alpha)}(\cF)) = \ch_{\mathrm{orig}}(\cF), \] 
where the third equality is Proposition~\ref{prop:etaatyiah}.
\end{proof}

\section{The main conjecture}
\label{sec:mainconj}

In this section we discuss the main conjecture and ways to approach
its proof.

\subsection{}
It was argued by Kontsevich~\cite{Kon} and Shoikhet~\cite{Sho} that
the isomorphisms arising from the Hochschild-Kostant-Rosenberg do {\em
not} respect the natural structures that exist on the Hochschild and
harmonic structures, respectively.  However, as a consequence of
Kontsevich's famous proof of the formality conjecture, he was able to
prove that correcting the $I^\HKR$ isomorphism by a factor of
$\Td_X^{1/2}$ yields a ring isomorphism:
\begin{theorem}[{\cite[Claim 8.4]{Kon}}]
\label{thm:kon}
Let $I^K$ be the composite isomorphism
\[
\begin{diagram}[height=2em,width=2em,labelstyle=\scriptstyle]
I^K:HH^*(X) & \rTo^{I^\HKR} & HT^*(X) & \rTo^{\vee \Td_X^{-1/2}} &
HT^*(X).
\end{diagram}
\]
Then $I^K$ is a ring isomorphism.
\end{theorem}

\subsection{}
Observe that the way the $I^\HKR$ isomorphism was defined, $I^K$ can
be defined with the same definition, but using a modified
Hochschild-Kostant-Rosenberg isomorphism
\[ I':\Delta^*\cO_\Delta \stackrel{\sim}{\lra} 
\bigoplus_i \Omega_X^i[i], \]
given by
\[ 
\begin{diagram}[height=2em,width=2em,labelstyle=\scriptstyle]
I':\Delta^*\cO_\Delta & \rTo^I & \bigoplus_i
\Omega_X^i[i] & \rTo^{\wedge\Td_X^{1/2}} & \bigoplus_i
\Omega_X^i[i].
\end{diagram}
\]
Here, by $\wedge \Td_X^{1/2}$ we have denoted the morphisms
\[
\begin{diagram}[height=2em,width=2em,labelstyle=\scriptstyle]
\Omega_X^j[j] & \rTo^{\Omega_X^j[j]\otimes \Td_X^{1/2}} & \bigoplus_i
\Omega_X^{i+j}[i+j],
\end{diagram}
\]
where 
\[ \Td_X^{1/2}:\cO_X \ra \bigoplus_i \Omega_X^i[i]\]
is the map that corresponds to 
\[ \Td_X^{1/2} \in \bigoplus_i H^i(X, \Omega_X^i) =
\Hom_X(\cO_X, \bigoplus_i \Omega_X^i[i]). \]

\subsection{}
\label{subsec:samemod}
The moral of Kontsevich's result is that $I$ is the ``wrong''
isomorphism to use, and the correct one is $I'$.  With this
replacement, $I_\HKR$ gets replaced by 
\[
\begin{diagram}[height=2em,width=2em,labelstyle=\scriptstyle]
I_K:HH_*(X) & \rTo^{I_\HKR} & H\Omega_*(X) & \rTo^{\wedge \Td_X^{1/2}}
& H\Omega_*(X).
\end{diagram}
\]
Not surprisingly, this matches well with the definition of the Mukai
vector: if we use $I$ and take Theorem~\ref{thm:chernhar} as our
definition of differential forms-valued Chern character, we get back
the classic definition of the Chern character; replacing $I$ by $I'$
replaces this classic Chern character with the Mukai vector
\[ v(\cF) = \ch(\cF) \wedge \Td_X^{1/2}, \]
which we saw in Sections~\ref{sec:prel} and~\ref{sec:harMukprod} is
better behaved from a functorial point of view.

\subsection{}
These observations, combined with the fact that all the properties of 
the Hochschild and the harmonic structures appear to match, lead us to
state the following conjecture:
\begin{conjecture}
\label{conj:mainconj}
The maps $(I^K, I_K)$ form an isomorphism between the
Hochschild and the harmonic structures of a compact smooth
space $X$.
\end{conjecture}

Observe that this conjecture includes, as a particular case,
Kontsevich's Theorem~\ref{thm:kon}.  An important consequence of this
conjecture will be discussed in~\cite{CalHH3}.

\subsection{}
We conclude with a remark on a possible approach to proving
Conjecture~\ref{conj:mainconj}.  For simplicity we restrict our
attention to a discussion of the isomorphism on cohomology (where we
know the conjecture is true by Kontsevich's result).  Consider the
sequence of morphisms
\[
\begin{diagram}[height=2em,width=3em,nohcheck,labelstyle=\scriptstyle,silent]
\Hom_X^*(\bigoplus \Omega_X^i[i], \bigoplus \Omega_X^i[i]) & \rTo^{I}
& \Hom_X^*(\Delta^*\cO_\Delta, \Delta^*\cO_\Delta) & & \\
\dTo^{p} & & \dTo^{-\circ \eta} & \luTo^{\Delta^*} & \\
\Hom_X^*(\bigoplus \Omega_X^i[i], \cO_X) & \rTo^I &
\Hom_{X}^*(\Delta^*\cO_\Delta, \cO_X) & \rTo^{\Delta_*(-)\circ 
\eta\ \ } & \Hom^*_{X\times X}(\cO_\Delta, \cO_\Delta) \\
\dEqual & & & & \dEqual \\
HT^*(X) & &\lTo^{I^\HKR} & & HH^*(X).
\end{diagram}
\]
The maps labeled $I$ are isomorphisms induced by $I$; the arrow
$\Delta_*(-)\circ \eta$ is the adjunction isomorphism.  The map $p$ is
the projection of a matrix in $\Hom_X^*(\bigoplus \Omega_X^i[i],
\bigoplus \Omega_X^i[i])$ onto its last column $\Hom_X^*(\bigoplus
\Omega_X^i[i], \cO_X)$.  (The convention that we use is that morphisms of
small degree appear at the {\em bottom} or {\em right} of column
vectors/matrices.)

Observe that all the vector spaces in the diagram have ring
structures, but only the top two and rightmost two have the ring
structure given by the Yoneda product.  Also, note that the arrows
between these rings are obviously ring homomorphisms.

We are interested in the map 
\[ e:\Hom_X^j(\bigoplus \Omega_X^i[i], \cO_X) \ra \Hom_X^j(\bigoplus
\Omega_X^i[i], \bigoplus \Omega_X^i[i]) \] 
which takes a column vector to a matrix, by the formula
\[
\left (
\begin{array}{c}
v_n\\
v_{n-1}\\
v_{n-2}\\
\vdots\\
v_0
\end{array}
\right ) \stackrel{e}{\mapsto} 
\left (
\begin{array}{ccccc}
v_0 & v_1 & v_2 & \cdots & v_n \\
0 & v_0 & v_1 & \cdots & v_{n-1} \\
0 & 0 & v_0 & \cdots & v_{n-2} \\
 & & \vdots & & \\
0 & 0 & 0 & \cdots & v_0
\end{array}
\right ).
\]
(For simplicity, at this point assume that we are only dealing with
{\em homogeneous} elements in $\Hom_X^*(\bigoplus \Omega_X^i[i],
\cO_X)$.)  It is easy to check that what we think of as
``multiplication'' in $\Hom_X^*(\bigoplus \Omega_X^i[i], \cO_X)$ is
the product
\[ v*v' = p(e(v)\circ e(v')). \]

There is another map $e'$ which takes a column vector and fills it up
to a square matrix $e'(v)$.  It is the map obtained by starting with
$v\in \Hom_X^j(\bigoplus \Omega_X^i[i], \cO_X)$ and following the
arrows around the diagram to get $e'(v)\in \Hom_X^j(\bigoplus
\Omega_X^i[i], \bigoplus \Omega_X^i[i])$.  The fact that $p\circ
e'$ is the identity means that the last column of $e'(v)$ is precisely
$v$.

To prove that $I^\HKR$ is a ring isomorphism, it would suffice to show
that $e'=e$.  Unfortunately, Kontsevich's argument shows that this is
not the case.  The same argument, however, shows that if we
repeat the above analysis with $I$ replaced by $I'$ (and $I^\HKR$
replaced by $I^K$) we do get a ring homomorphism.  This leads us to
state the following conjecture:
\begin{conjecture}
Replacing $I$ by $I'$ in the above analysis yields $e=e'$.
\end{conjecture}

A proof of this conjecture, apart from providing a different proof of
Kontsevich's result, would likely generalize to a proof of
Conjecture~\ref{conj:mainconj}.

\bigskip \noindent
\small\textsc{Department of Mathematics, \\
University of Philadelphia, \\
Philadelphia, PA 19104-6395, USA} \\
{\em e-mail: }{\tt andreic@math.upenn.edu}

\end{document}